\documentclass[12pt]{amsart}

\usepackage{amsgen,amsmath,amstext,amsbsy,amsopn,amsfonts,amssymb}

\def\char{{\rm char}}
\def\Proof{\noindent{\sl Proof.}\ }
\def\qed{{\hfill $\Box$ \medbreak}}

\newtheorem{defi}{Definition}[section]
\newtheorem{lem}[defi]{Lemma}
\newtheorem{thm}[defi]{Theorem}
\newtheorem{cor}[defi]{Corollary}
\newtheorem{remark}[defi]{Remark}
\newtheorem{prop}[defi]{Proposition}

\DeclareMathOperator{\la}{\langle}
\DeclareMathOperator{\ra}{\rangle}
\DeclareMathOperator{\su}{\subseteq}

\DeclareMathOperator{\ad}{ad}

\DeclareMathOperator{\F}{\mathbb{F}}

\begin{document}

\title[]{Lie identities on symmetric elements of restricted enveloping algebras }
\author{\textsc{Salvatore Siciliano}}
\address{Dipartimento di Matematica ``E. De Giorgi", Universit\`{a} del Salento,
Via Provinciale Lecce--Arnesano, 73100--Lecce, Italy}
\email{salvatore.siciliano@unisalento.it}

\author{\textsc{Hamid Usefi}}
\address{Department of Mathematics, University of Toronto, Bahen Centre, 40 St. George St.,
Toronto, Ontario, Canada, M5S 2E4}
\email{usefi@math.toronto.edu}

\begin{abstract}
Let $L$ be a restricted Lie algebra over a field of characteristic $p>2$ and denote by $u(L)$ its restricted enveloping algebra. 
We determine the conditions under which the set of symmetric elements of $u(L)$ with respect to the principal involution is Lie solvable, Lie nilpotent, or bounded Lie Engel. 

\end{abstract}
\subjclass[2010]{16S30,  16W10, 17B50, 17B60}
\date{\today}

\maketitle{}

\section{Introduction}

Let $A$ be an algebra with involution $\ast$ over a field $\F$. We denote by $A^+:=\{x\in A\vert\,x^\ast=x\}$ the set of symmetric elements of $A$ under $\ast$ and by $A^-:=\{x\in A\vert\,x^\ast=-x\}$ the set of skew-symmetric elements.  A  question of general interest is which properties of $A^+$ or $A^-$ can be lifted to the whole algebra $A$. The history of this problem goes back to Herstein \cite{H67, H76} where he had conjectured that if  the symmetric or skew-symmetric elements of  a ring $R$ satisfy a polynomial identity, then so does $R$. Notably this conjecture was proved  by  Amitsur in  \cite{A68} and subsequently generalized by himself in  \cite{A69}.

Now consider the group algebra $\F G$ of a group $G$ under the canonical involution induced by $g \mapsto g^{-1}$, for every $g\in G$. There has been an intensive investigation devoted to demonstrate the extent to which the symmetric or skew-symmetric elements of  $\F G$ under the canonical involution  determine the structure of the group algebra and there has been special attention on Lie identities. In particular, the chacterization of groups $G$ for which 
$\F G^-$ or $\F G^+$ is Lie nilpotent was carried out by Giambruno,  Sehgal, and Lee in \cite{GS93, GS06, L99}. Furthermore, if either $\F G^-$ or $\F G^+$ is bounded Lie Engel, and $G$ is devoid of $2$-elements, Lee in \cite{ L00} showed that $\F G$ is bounded Lie Engel. He also classified remaining groups for which $\F G^+$ is bounded Lie Engel. 
The Lie solvable case was considered in \cite{LSS, LSS1}, however a complete answer to this case seems still under way.

Now, let $L$ be a restricted Lie algebra  over a field ${\mathbb F}$ of characteristic $p>2$ and let $u(L)$ be the restricted enveloping algebra of $L$.  
We denote by $\top$  the \emph{principal involution} of $u(L)$, that is, the unique ${\mathbb F}$-antiautomorphism of $u(L)$ such that $x^\top=-x$ for every $x$ in $L$. We recall that $\top$ is just the antipode of the  ${\mathbb F}$-Hopf algebra $u(L)$.

Recently, the first author in \cite{Sal}  established the conditions under which  $u(L)^-$ is Lie solvable, Lie nilpotent or  bounded Lie Engel. 
 
In this paper we consider the symmetric case. Unlike the skew-symmetric case it is not clear a priori that, for example, when $u(L)^+$ is Lie nilpotent then $L$ is nilpotent. 
The symmetric elements do not form a Lie subalgebra of $u(L)$  in general (but they form a Jordan subalgebra under the Jordan bracket $x\circ y=\frac{1}{2} (xy+ yx)$).
However, despite the group ring case we present  a complete answer and yet  the  proofs are different and more involved than  the skew-symmetric case. Before stating the main results we recall the following definitions.
An element $x$ of $L$ is $p$-nilpotent if $x^{[p]^m}=0$ for some $m\geq 1$; a subset $S$ of $L$ is  $p$-nilpotent if $S^{[p]^m}=\{x^{[p]^m}\vert\, x\in S \}=0$ for some $m\geq 1$. The derived subalgebra of $L$ is denoted by $L^\prime$. 

\begin{thm}\label{bengel}
Let $L$ be a restricted Lie algebra over a field ${\mathbb F}$ of characteristic $p>2$. Then the following conditions are equivalent:
\begin{enumerate}
\item [{\normalfont 1)}] $u(L)^+$ is bounded Lie Engel;
\item [{\normalfont 2)}] $u(L)$  is bounded Lie Engel;
\item [{\normalfont 3)}] $L$ is nilpotent, $L^\prime$ is $p$-nilpotent, and $L$
contains a restricted ideal $I$ such that $L/I$ and $I^\prime$ are finite-dimensional.
\end{enumerate}
\end{thm}

\begin{thm}\label{nilpotent}
Let $L$ be a restricted Lie algebra over a field ${\mathbb F}$ of characteristic $p>2$. Then the following conditions are equivalent:
\begin{enumerate}
\item [{\normalfont 1)}] $u(L)^+$ is Lie nilpotent;
\item [{\normalfont 2)}] $u(L)$ is Lie nilpotent;
\item [{\normalfont 3)}] $L$ is nilpotent and $L^\prime$ is finite-dimensional and $p$-nilpotent.
\end{enumerate}
\end{thm}

\begin{thm}\label{solvable}
Let $L$ be a restricted Lie algebra over a field ${\mathbb F}$ of characteristic $p>2$. Then the following conditions are equivalent:
\begin{enumerate}
\item [{\normalfont 1)}] $u(L)^+$ is Lie solvable;
\item [{\normalfont 2)}] $u(L)$ is Lie solvable;
\item [{\normalfont 3)}] $L^\prime$ is finite-dimensional and $p$-nilpotent.
\end{enumerate}

\end{thm}

 The equivalence of 2) and 3) in Theorems 1.1-1.3 was established in \cite{RS}. Our main contribution is to prove that 2) implies 3) in these theorems. In combination with \cite{Sal}, we can  conclude that, in odd characteristic, if either $u(L)^+$ or $u(L)^-$ is Lie solvable (respectively, bounded Lie Engel or Lie nilpotent) then so is the whole algebra $u(L)$. Such conclusions are no longer true in characteristic $2$.

Finally, let $L$ be an arbitrary Lie algebra over a field of characteristic different from $2$ and denote by $U(L)$ the ordinary enveloping algebra of $L$. A further consequence of our main results is that  $U(L)^+$  (under the principal involution) is Lie solvable or bounded Lie Engel only when $L$ is abelian.

\section{Preliminaries}

Throughout the paper  all restricted Lie algebras are defined over a  field $\F$ of characteristic $p>2$.

Let $A$ be an associative algebra over $\F$ with an involution. 
The Lie bracket on $A$ is defined by $[x,y]=xy-yx$, for every $x,y\in A$. Longer Lie commutators of $A$ are recursively defined as follows: $[x_1,x_2,\ldots,x_{n+1}]=[[x_1,x_2,\ldots,x_n],x_{n+1}]$
and $[x,_{n+1} y]=[[x,_ny],y]$.  A subset $S$ of $A$ is said to be Lie nilpotent if there exists a positive integer $n$ such that 
$$[x_1,\ldots,x_n]=0$$
 for every $x_1,\ldots,x_n\in S$, while  $S$ is said to be bounded Lie Engel if there is an $n$ such that 
 $$[x,_ny]=0$$
  for every $x,y\in S$. Also, following \cite{L10}, we put $[x_1,x_2]^o =[x_1,x_2]$ and 
$$
[x_1,x_2,\ldots,x_{2^{n+1}}]^o=[[x_1,\ldots,x_{2^n}]^o,
[x_{2^n+1},\ldots,x_{2^{n+1}}]^o].
$$ 
The subspace spanned by all
 $[a_1,\ldots,a_{2^n}]^o$, where $a_i\in S$, is denoted by $\delta_n(S)$.
The subset $S$  is said to be Lie solvable if there exists an $n$ such that 
$\delta_n(S)=0$.

Moreover, we denote by $\gamma_n(L)$ ($n\geq 1$) and $\delta_n(L)$ ($n\geq 0$) the terms of the lower central series and derived series of $L$, respectively. Also, for a subset $S$ of $L$ we denote by $S_p$ the restricted subalgebra generated by $S$. 
A polynomial of the form $x^{ p^k}+\alpha_1 x^{ p^{k-1}}+\cdots+\alpha_k x\in \F[x]$ is called a $p$-polynomial. An element $x$ of $L$ is said to be $p$-algebraic if $\dim_{\F} \la x \ra_p<\infty$. Finally, $Z(L)$ is the center of $L$.

\begin{thm}[\cite{A68}]\label{amitsur} 
Let $A$ be an associative algebra with involution. If $A^+$ (or $A^-$) satisfies a 
polynomial identity then $A$ satisfies a polynomial identity.
\end{thm}

Restricted Lie algebras whose enveloping algebras satisfy polynomial identities are
characterized by Passman \cite{Pa} and Petrogradski \cite{Pe}.
A theorem of Kukin  (see \cite[\S 4, Lemma 1.10]{K}) states  that 
any  restricted subalgebra  $H$ of $L$ of finite codimension contains 
an ideal of finite codimension in $L$.
Combining these results we get the following:

\begin{thm}\label{passman} 
Let $L$ be a restricted Lie algebra over a field of characteristic $p>0$.
Then $u(L)$ satisfies a polynomial identity if and only if $L$ has a restricted  ideal 
$A$ of finite codimension in $L$ such that $A$ is nilpotent of class two and $A^\prime$ is finite-dimensional and $p$-nilpotent.
\end{thm}

\begin{thm}[\cite{St}]\label{stewart}
Let $L$ be a Lie algebra with a nilpotent ideal $M$ such that
$L/M^\prime$ is nilpotent. Then $L$ is nilpotent.
\end{thm}

\begin{cor}\label{stewart-cor}
Let $L$ be a restricted Lie algebra with a nilpotent ideal $M$ such that
$L/(M^\prime)_p$ is nilpotent. Then $L$ is nilpotent.
\end{cor}
\Proof Suppose that $\gamma_{c+1}(L)\su (M^\prime)_p$. Then $\gamma_{c+2}(L)\su [ (M^\prime)_p, L]\su M^\prime$.
Thus, $L/M^\prime$ is a nilpotent Lie algebra and it follows from Theorem \ref{stewart} that $L$ is nilpotent.
\qed

\begin{lem}\label{infinitecenter-lem}
Let $L$ be a restricted Lie algebra  such that  $u(L)^+$ is Lie solvable.
If $Z(L)$ is infinite-dimensional,   then $u(L)$ is Lie solvable.
\end{lem}
\Proof The proofs of Lemma 1 and Corollary 1 in \cite{Sal} work in the symmetric situation as well (by using $x_i+x_i^\ast$ instead of $x_i-x_i^\ast$). \qed

A similar conclusion also holds for the  Lie nilpotence case. Now suppose that $u(L)^+$ is Lie solvable
and $L$ is nilpotent of class $c$, say.  If $\gamma_c(L)_p$ is infinite-dimensional, then because  $\gamma_c(L)_p$ is a central ideal of $L$, we deduce by Lemma \ref{infinitecenter-lem}  that
$u(L)$ itself is Lie solvable.  But   then  $L^\prime$ is  finite-dimensional and $p$-nilpotent, see  \cite[Theorem 1.3]{RS} or 
\cite[Lemma 4.4]{U}. This contradiction yields the  following:
\begin{cor}\label{infinitecenter}
Let $L$ be a nilpotent restricted Lie algebra  of class $c\geq 2$. If  $u(L)^+$ is Lie solvable, then $\gamma_c(L)_p$ is finite-dimensional.
\end{cor}

The following can be  deduced from classical results. For the proof see, for example,   \cite{RS} or \cite{U}.

\begin{lem}\label{u(L)-Lie-solvable}
If $u(L)$ is Lie solvable then $L^\prime$ is $p$-nilpotent.
\end{lem}

In the sequel, we shall freely use the following fact.

\begin{remark} \emph{ Let $A$ be an algebra over a field of characteristic not $2$ with involution $\ast$. If $J$ is a $\ast$-invariant ideal of $A$ then  
$\big( A/J \big)^+$ with respect to the induced involution coincides with the image of $A^+$ under the canonical map $A \to A/J$. In particular, this applies to the case where $A=u(L)$ and $J$ is the associative ideal generated by a restricted ideal of $L$.}
\end{remark}
\section{Proofs of the main results}

\begin{lem}\label{f.d.-nilpotent}
Let $L=\langle H \rangle_p$ where $H$ is a finite-dimensional Lie algebra. If $u(L)^{+}$ is bounded Lie Engel then $L$ is nilpotent.
\end{lem}
\Proof Clearly, it is enough to show that $H$ is nilpotent.  By Engel's Theorem it suffices to prove that $\ad x$ is nilpotent, for every $x\in H$. 
Suppose to the contrary that $\ad x$ is not nilpotent, for some $x\in H$. 
Let $\bar \F$ be the algebraic closure of $\F$. Since $\ad (x\otimes 1)$ is not a  nilpotent operator, 
there exists a nonzero $ b\in H\bigotimes \bar \F$ such that $[ x\otimes 1 ,  b]=\alpha  b$, where  $\alpha$ is a non-zero eigenvalue of $\ad (x\otimes 1)$. We may replace $x$ with $\alpha^{-1} x$ to assume that $[ x\otimes 1 ,  b]= b$.
Let $b=\sum_{k=1}^t b_k\otimes \alpha_k$, where $b_k\in H$ and $\alpha_k\in \F$. Then 
$$
b^2=\sum_{1\leq i< j\leq t} (b_ib_j+b_jb_i)\otimes \alpha_i\alpha_j +\sum_{k=1}^t b_k^2\otimes \alpha_k^2.
$$
Note that $b_ib_j+b_jb_i, x^2\in u(L)^+$. There exists an integer $m$ such   that
$$
[b_ib_j+b_jb_i, _{p^m} x^2]=0
$$
for all $i,j$. Similarly, $[b_k^2, _{p^m} x^2]=0$, for every $k$. It follows that 
\begin{align*}
[ b^2, _{p^m} x^2\otimes 1]=& \sum_{1\leq i< j\leq t} [(b_ib_j+b_jb_i), _{p^m} x^2]\otimes \alpha_i\alpha_j +\sum_{k=1}^t [b_k^2, _{p^m} x^2]\otimes \alpha_k^2=0.
\end{align*}
On the other hand, we have 
$$[x^2\otimes 1,  b^2]=[(x\otimes 1)^2,  b^2]=4 b^2((x+1)\otimes 1).
$$
Thus,
$$
[ b^2, _{p^m} \frac{x^2\otimes 1}{4}]=-  b^2((x+1)^{p^m}\otimes 1)=-  b^2(x\otimes 1)^{p^m}- b^2=0
$$
By the PBW Theorem for restricted Lie algebras (see e.g. \cite[\S 2, Theorem 5.1]{SF}) this is possible only if $x$ is $p$-nilpotent which implies that
$\ad x$ is nilpotent, a contradiction.
\qed

\begin{lem}\label{metabelian-nilpotent}
Let $L$ be a metabelian restricted Lie algebra containing an abelian ideal $N$ of finite codimension. If  $u(L)^{+}$ is bounded Lie Engel,  then $L$ is nilpotent.
\end{lem}
\Proof We argue by induction on $\dim L/N$. 
Note that $L/N$ is nilpotent, by Lemma \ref{f.d.-nilpotent}.
Let $M=N+(L^\prime)_p$.
Then $M$ is metabelian
and since $L/N$ is nilpotent, we have $\dim M/N<\dim L/N$.
Thus, by induction on $\dim L/N$, we can assume that $M$ is nilpotent.
In order to prove $L$ is nilpotent,  it is enough, by Corollary \ref{stewart-cor},
 to show that $L/(M^\prime)_p$ is nilpotent.
So we replace $L$ with $L/(M^\prime)_p$ to assume that $[L^\prime, N]=0$. We also replace $N$ with $N+(L^\prime)_p$ to assume that $L^\prime\su N$.
Now let $x\in L\setminus N$.
Since $L/N$ is finite-dimensional, there exists a $p$-polynomial $f$ such that
$f(x)\in N$. It is easy to see that, for every $a\in N$, the Lie algebra $H$ generated by $x$ and $a$ is of dimension at most $\deg(f)+1$.
So, by Lemma  \ref{f.d.-nilpotent},   $H_p$ is nilpotent. Indeed, the nilpotency class of $H_p$ depends only on $\dim H$.
We deduce that $[N,_m x]=0$, where $m$ depends only on $x$. Hence the restricted subalgebra $I$ generated by $N$ and $x$ is nilpotent.
Note that $I$ is in fact a restricted ideal of $L$. By Corollary \ref{stewart-cor}, we can replace $L$ with $L/(I^\prime)_p$.
But then $\dim L/I< \dim L/N$. So, by the induction hypothesis, $L$ is nilpotent, as required.
\qed

\begin{prop}\label{L-nilpotent}
Let $L$ be a  restricted Lie algebra  such that $u(L)^{+}$ is bounded Lie Engel. Then $L$ is nilpotent.
\end{prop}
\Proof
It follows from Theorem \ref{amitsur} that $u(L)$ is PI. Thus, by Theorem \ref{passman}, there exists a nilpotent restricted ideal $N$ of $L$ of finite codimension. Hence, by Corollary \ref{stewart-cor}, it is enough to prove that $L/(N^\prime)_p$ is nilpotent.
Thus, we can assume that $L$ has an abelian ideal $N$ of finite codimension.
Since $L/N$ is finite-dimensional, by Lemma \ref{f.d.-nilpotent} we see that $L/N$ is nilpotent.
Thus $L$ is solvable.
 Let $d$ be the derived length of $L$ and consider $M=\delta_{d-2}(L)_p$.
By Lemma \ref{metabelian-nilpotent}, we know that $M$ is nilpotent. Moreover, by induction on the derived length we have that $L/(M^\prime)_p=L/\delta_{d-1}(L)_p$ is nilpotent. Therefore, by  Corollary \ref{stewart-cor}, we conclude that $L$ is nilpotent.
\qed

\begin{cor}
Let $L$ be a  restricted Lie algebra. The symmetric elements of $u(L)$ commute if and only if $L$ is abelian.
\end{cor}
\Proof The sufficiency is obvious. Suppose $u(L)^{+}$ is commutative and that, if possible,  $L$ is not abelian.
Since, by Proposition \ref{L-nilpotent}, $L$ is nilpotent, we can find 
non-commuting elements $x, y\in L$ such that $z=[x, y]$ is central.
But then $x^2$ and $2xy-z$ are symmetric elements of $u(L)$ and yet  $[x^2,2xy-z]=4x^2z$ is not zero, a contradiction.
\qed

\begin{prop}\label{derivedpnilp}
Let $L$ be a  restricted Lie algebra such that $u(L)^{+}$ is bounded Lie Engel. Then $L^\prime$ is $p$-nilpotent.
\end{prop}
\Proof
By hypothesis, there exists an $m$ such that 
 $[a,_mb]=0$ for every $a,b\in u(L)^{+}$. Suppose, if possible, that  $L^\prime$ is not $p$-nilpotent.
By Proposition \ref{L-nilpotent} we know $L$ is nilpotent, hence there exists a minimal
$n> 2$ such that $\gamma_n(L)_p$ is $p$-nilpotent.
 Put ${\mathfrak L}:=L/\gamma_n(L)_p$.
Let $t$ be an integer such that $p^t\geq m$. Then from the assumption it follows
that there exist  $x,y\in {\mathfrak L}$ such that $z=[x,y]\in Z({\mathfrak L})$ and $z^{[p]^t}\neq 0$. Clearly, the element $x^2$ and $2xy-z$ of $u({\mathfrak L})$ are symmetric.
Moreover, it can be seen by an easy induction that for every $r>0$ one has
\begin{equation}\label{engel}
[x^2,_{r}2xy-z]=4^rx^2z^{r}.
\end{equation}
Since $p^t\geq m$ and $p\neq 2$, from relation (\ref{engel}) we deduce that $x^2z^{p^t}=0$. On the other hand, since $x$ and $z^{[p]^t}$ are  ${\mathbb F}$-linearly independent (as $z^{[p]^t}$ is central whereas $x$ is not), the last conclusion contradicts the PBW Theorem, completing the proof.
\qed

It is now a simple matter to prove Theorems \ref{bengel} and \ref{nilpotent}.
\smallskip

\noindent{\sl Proof of Theorem \ref{bengel}.}
Of course, in view of Theorem 2 of \cite{RS}, it is enough to show that 1) implies 3).
By Lemmas \ref{L-nilpotent} and \ref{derivedpnilp}, we already know that $L$ is nilpotent and $L^\prime$ is $p$-nilpotent.
 Moreover, by Theorem \ref{amitsur}, $u(L)$ satisfies a polynomial identity and then,  by Theorem \ref{passman}, $L$ contains a restricted ideal $I$ such that $L/I$ and $I^\prime$ are finite-dimensional. This finishes the proof.
\qed

\noindent{\sl Proof of Theorem \ref{nilpotent}.}
It is enough, by Theorem 1 of \cite{RS}, to show that $1)$ implies $3)$. So, suppose $u(L)^+$ is Lie nilpotent. Then, by Theorem  \ref{bengel}, we just need to show that $L^\prime$ is finite-dimensional. By Theorem \ref{bengel}, we know that $L$ is nilpotent of class $c$, say. Now, by Corollary \ref{infinitecenter},  $\gamma_c(L)_p$ is  finite-dimensional. We can proceed by induction on $c$ to conclude that $L^\prime+\gamma_c(L)_p/\gamma_c(L)_p$ is finite-dimensional, hence $L^\prime$ is  finite-dimensional, too. \qed   

\begin{prop}\label{L-solvable}
Let $L$ be a  restricted Lie algebra  such that $u(L)^{+}$ is Lie solvable. Then $L$ is solvable.
\end{prop}
\Proof We can suppose, without loss of generality, that the ground field $\F$ is algebraically closed.
Moreover, by Theorem \ref{amitsur}, $u(L)$ is PI and then, by Theorem \ref{passman}, there exists a nilpotent restricted ideal $N$ of $L$ of finite codimension. Thus we can replace $L$ by $L/N$ and assume that $L$ is finite-dimensional. Suppose that $L$ is a counterexample of minimal dimension. Now, if $M$ is a non-trivial restricted ideal of $L$, then $u(J)^+$ and $u(L/M)^+$ are Lie solvable. Hence the minimality of $L$ forces that both $M$ and $L/M$ are solvable, so that $L$ is solvable, a contradiction. Consequently, $L$ has no non-trivial restricted ideal.  Now, we claim that $L^\prime$ is a simple Lie algebra. Let $I$ be a nonzero ideal of $L^\prime$. As $(L^\prime)_p$ is a nonzero restricted ideal of $L$ we have  $(L^\prime)_p=L$ and then $I$ is an ideal of $L$.
It follows that $I_p=L$ and then $L^\prime=[I_p,I_p]=[I,I]\subseteq I$, so that $I=L^\prime$, as claimed. In particular, $L^\prime$ is a simple restricted $L$-module.
Denote by ${\mathcal J}$ the Jacobson radical of $u(L)$. Note that the simple modules of $u(L)$ and $u(L)/{\mathcal J}$ are the same. Since ${\mathcal J}$ is $\top$-invariant we have that $\left( u(L)/{\mathcal{J}} \right)^+$ under the induced involution is Lie solvable. As a consequence, if $u(L)/{\mathcal J}\cong \oplus_{i=1}^r M_{n_i}(\F)$ is the the Wedderburn decomposition of $u(L)/{\mathcal J}$, in view of Lemma $3.2$ of \cite{LSS} we have that $n_i\leq 2$ for every $i$. Therefore every irreducible representation of $u(L)$ has degree $1$ or $2$. In particular $L^\prime$ has dimension less than $3$, so that $L$ is solvable, a contradiction.
\qed

\begin{lem}\label{solvable-p-nil}
Let $L$ be a finite-dimensional restricted Lie algebra  such that $u(L)^{+}$ is Lie solvable. Then $L^\prime$ is $p$-nilpotent.
\end{lem}
\Proof By Proposition \ref{L-solvable}, we know that $L$ is solvable.
By induction on the derived length of $L$, we can assume $L$ is metabelian. Furthermore, we can suppose that the ground field $\F$ is algebraically closed.
Let ${\mathcal J}$ be the Jacobson radical of $u(L)$ and consider the restricted Lie algebra $H=L+{\mathcal J}/{\mathcal J}$.
Since ${\mathcal J}$ is a nilpotent ideal of $u(L)$, it suffices to show that $H^\prime$ is $p$-nilpotent.
Thus, we may assume that $H$ is not abelian.
Note that $u(L)/{\mathcal J}\cong \F\oplus\cdots \F\oplus M_2(\F)\oplus\cdots \oplus M_2(\F)$, by \cite{LSS}.
Without loss of generality we can assume that $H$ is a restricted Lie subalgebra of $M_2(\F)$.
Since $M_2(\F)$ is not Lie metabelian, we have $\dim H\leq 3$. If $\dim H=2$, then
there exists $x, y\in H$ such that  with $[x,y]=x$. In this case
we must have $x^p=0$ and $y^p=y$, so we are done.
Assume then $\dim H=3$. It is clear that the identity matrix $I_2$ must be in $H$,
otherwise $M_2(\F)$ is spanned by $H$ and $I_2$ so that $M_2(\F)$ is Lie metabelian, which is not possible. Note
that $H^\prime \subseteq [M_2(\F),M_2(\F)]=\text{sl}_2(\F)$. Since $H^\prime$ is abelian
and $\text{sl}_2(\F)$ has no 2-dimensional abelian subalgebra, it follows that $\dim H^\prime=1$.
Let $x\in H^\prime\su \text{sl}_2(\F)$. Since $\text{tr}(x)=0$, it follows that $x$ and $I_2$ are linearly independent.
Let $y\in H$ so that $x, y$ and $I_2$ span $H$. Since $H$ is not abelian, we can assume
 $[x,y]=x$. But then   $x^p$ is a central element in $H$ and so $x^p=\alpha I_2$, for some $\alpha\in \F$. Since $x\in \text{sl}_2(\F)$ and $p\geq 3$,
$\text{tr}(x^p)=0$. Thus, $\alpha=0$ and so $x^p=0$.
\qed

\begin{lem}\label{nilpotent-p-nil}
Let $L$ be a nilpotent restricted Lie algebra   such that $u(L)^{+}$ is Lie solvable. Then $L^\prime$ is finite-dimensional and $p$-nilpotent.
\end{lem}
\Proof Let $c$ be the nilpotence class of $L$.
First we prove that $\gamma_c(L)$ is $p$-nil.
So let $x\in L$ and $y\in \gamma_{c-1}(L)$  and let $H$ be the restricted 
subalgebra of $L$ generated by $x$ and $y$.
 Let ${\mathcal J}$ be the Jacobson radical of $u(H)$.
By the Razmyslov-Kemer-Braun Theorem (see e.g. \cite{B}) ${\mathcal J}$ is a nilpotent ideal of $u(H)$.
 Note that $x^2, y^2$ and $u=2xy-z$ are symmetric.
In $u(H)$ we have $[x, 2xy-z]=2xz$. So,  $[x^2, u]=4x^2z$.
We can similarly deduce that,  $[y^2, u]=-4y^2z$. Thus,
\begin{align*} 
[[x^2, u], [y^2, u]]=16[x^2z, y^2z]=-64xyz^3+32z^4.
\end{align*}
By Proposition 2.6 of  \cite{LSS}, $(u(H)/{\mathcal J})^+$ is Lie metabelian. It follows that $2xyz^3-z^4\in {\mathcal J}$.
Note that 
$$
[[x, 2xyz^3-z^4], y]=2z^5\in {\mathcal J}.
$$
Thus, there exists an integer $k$ such that 
$$
(z^{p^k})^5=0, 
$$
which implies that $z$ must be $p$-nilpotent.
Note that $\gamma_c(L)_p$  is finite-dimensional, by  Corollary \ref{infinitecenter}.
Now we consider $L/\gamma_c(L)_p$. By induction on $c$,  we have that $L^\prime+\gamma_c(L)_p/\gamma_c(L)_p$ is finite-dimensional and $p$-nilpotent. Combining with the fact that
$\gamma_c(L)_p$ is also finite-dimensional and $p$-nilpotent yields the required result.
\qed


\begin{lem}\label{fixedpoints}
 Let $L$ be a metabelian restricted Lie algebra such that $u(L)^+$ is Lie solvable.
 Then the  space  of fixed points of the action of $\ad x$ on $L^\prime$ is finite-dimensional, for every $x\in L$.
\end{lem}
\Proof Let $V$ be the subspace of $L^\prime$ consisting of the fixed points of the action of $\ad x$ on $L^\prime$.
Note that for every $r\geq 1$ and $a_1, \ldots, a_r\in V$, we have
 $$
 [x, a_1\cdots a_r]=ra_1\cdots a_r.
 $$
Thus, for every odd number $k$ and $a_1, \ldots, a_k\in V$, we have
\begin{align}\label{n=0}
2xa_1\cdots a_{k}-ka_1\cdots a_{k}=xa_1\cdots a_{k} +(xa_1\cdots a_{k})^\top \in u(L)^+.
\end{align}
Let $n$ be a positive integer and let 
  $m_1=2^{n+1}+2^n-2$, $m_2=2^{n+1}+2^n$ and $m_3=2^{n+1}+2^n+2$.
By the PBW Theorem, it is enough to show that for every $a_1,\ldots ,a_{m_i}\in V$ there exists $\alpha\in \F$ such that 
$$
xa_1a_2\cdots a_{m_i} + \alpha a_1a_2\cdots a_{m_i}\in \delta_n(u(L)^+),
$$
for $i=1,2,3$.
Using Equation (\ref{n=0}), the claim  can be easily checked for $n=0$.
Now we assume that $n\geq 1$ and we prove the induction step for $m_1$. So set $r=2^{n}+2^{n-1}-2$, $s=2^n+2^{n-1}$ and
let $a_1,\ldots, a_{r}, b_1,\ldots, b_{s}\in V$. By the induction hypothesis we have
\begin{align*}
u&=xa_1\cdots a_{r}+\alpha a_1\cdots a_{r}\in  \delta_{n}(u(L)^+),\\
v&=xb_1\cdots b_{s}+\beta b_1\cdots b_{s} \in  \delta_{n}(u(L)^+),
\end{align*}
for some $\alpha, \beta\in \F$. Thus,
\begin{align*}
[u, v]
=&2xa_1\cdots a_{r}b_1\cdots b_{s}+ (\beta s-\alpha r) a_1\cdots a_{r}  b_1\cdots b_{s}\in  \delta_{n+1}(u(L)^+).
\end{align*}
Similarly, for $m_2$ we take $r=2^n+2^{n-1}-2$ and $s=2^n+2^{n-1} +2$ and for
$m_3$ we take $r=2^n+2^{n-1}$ and $s=2^n+2^{n-1}+2$
and argue the same way to establish the inductive step.
\qed

\begin{lem}\label{codim-1}
 Let $L=\la y\ra_p+N$ where $N$ is an  abelian restricted ideal of finite codimension in $L$.
 If $u(L)^+$ is Lie solvable, then $L^\prime$ is finite-dimensional.
\end{lem}
\Proof
Without loss of generality we can suppose the ground field $\F$ is algebraically closed.
We argue by induction on $\dim L/N$. 
By \cite[\S 2, Theorem 3.6] {SF},  there exists a $p$-polynomial $g$ such that
$x=g(y)$ lies in $L\backslash N$ and satisfies either $x^{[p]}=0$ (mod $N$) or $x^{[p]}=x$ (mod $N$).
Thus,  $(\ad x)^p=0$  or $(\ad x)^p=\ad x$ in $L$.
Let $H=\la x, N\ra_p$. We claim that $(H^\prime)_p$ is finite-dimensional.
 If $(\ad x)^p=0$, then $H$ is nilpotent and the claim follows from  Lemma \ref{nilpotent-p-nil}. On the other hand, if $(\ad x)^p=\ad x$, since $H^\prime=[x,H]$ then $(\ad x)^{p-1}$ acts the identity on $H^\prime$. Therefore all the eigenvalues of  $\ad x$ on $H^\prime$ are  $(p-1)^{\text {th}}$ root of unity and $H^\prime$ decomposes as the sum of the relative eigenspaces. Now, let $\lambda$ be  an eigenvalue of $\ad x$ and   put $\bar{x}=\lambda^{-1}x$.
Then, by Lemma \ref{fixedpoints}, the vector space of fixed points of $\ad \bar{x}$ on $H^\prime$ is finite-dimensional, which is equivalent to saying that the eigenspace relative to the eigenvalue $\lambda$ of $\ad x$ is finite-dimensional. Hence,  $H^\prime$ is finite-dimensional. If there exists an element $u\in H^\prime$ such that $\la u\ra_p$ is not finite-dimensional, then the center of $H$ would be infinite-dimentional. This implies, by Lemma \ref{infinitecenter}, that $u(H)$ is Lie solvable and in turn, by Lemma \ref{u(L)-Lie-solvable}, $H^\prime$ is $p$-nilpotent, a contradiction.
We conclude that $(H^\prime)_p$ is finite-dimensional, as claimed. Now we replace $L$ with $L/(H^\prime)_p$ to assume that $H$ is abelian. 
But then $L$ has an  abelian restricted ideal $H$ of smaller codimension than $N$. Therefore it follows by induction that $L^\prime$ is finite-dimensional.
\qed

\noindent{\sl Proof of Theorem \ref{solvable}.} 
By Theorem 3 of \cite{RS}, it is enough to show that $1)$ implies $3)$.
Note that, by Proposition \ref{L-solvable},  $L$ is solvable. Hence, by induction on the derived length, we can assume that $L$ is metabelian.
By Theorem \ref{amitsur}, $u(L)$ is PI. Thus, by Theorem \ref{passman}, there exists a nilpotent restricted ideal $N$ of $L$ of finite codimension.
By Lemma \ref{nilpotent-p-nil}, $N^\prime$ is finite-dimensional and $p$-nilpotent.
So we can replace $L$ with $L/(N^\prime)_p$ to assume that $N$ is abelian. 
Now we argue by induction on $\dim L/N$.
Let $I$ be the restricted ideal generated by $L^\prime+N$.
Since $L$ is metabelian, it follows that $I$ is nilpotent of class two. Thus again, by Lemma \ref{nilpotent-p-nil}, we can replace $L$ with
$L/(I^\prime)_p$ to assume that $[N, L^\prime]=0$. 
Now let $x\in L\setminus N$
and denote by $H$ the restricted subalgebra generated by $x$ and $L^\prime+N$.
 By Lemma \ref{codim-1}, $H^\prime$ is finite-dimensional.
 We claim that $H^\prime$ is $p$-nilpotent.
 If $L^\prime+N$ has an element which  is not $p$-algebraic then  $Z(H)$ is infinite-dimensional and it follows  from Lemma \ref{infinitecenter} that
$u(H)$ is Lie solvable which in turn implies, by Lemma \ref{u(L)-Lie-solvable}, that $H^\prime$ is $p$-nilpotent. Hence we can assume that every element of $L^\prime+N$ is $p$-algebraic. Since $\dim L/N$ is finite, there exists a $p$-polynomial such that $z=f(x)\in N$.
Thus $z$ is $p$-algebraic and so is $x$. Now let
$y\in L^\prime+N$. Then note that
$M=\la x, y\ra_p$ is finite-dimensional. Hence, by Lemma \ref{solvable-p-nil}, 
 $M^\prime$ is $p$-nilpotent.
 We just showed that $[x,y]$ is $p$-nilpotent for every $y\in L^\prime+N$.
 Thus, $H^\prime$ is $p$-nilpotent. 
Now we can replace $L$ with $L/(H^\prime)_p$. 
So, $H$ is an abelian restricted ideal of $L$ and $\dim L/H < \dim L/N$.
Thus, by induction hypothesis,  $L^\prime$ is finite-dimensional and $p$-nilpotent, as required
\qed

\section{Concluding remarks}

Note that if $A$ is any algebra with involution then $[A^+, A^+]\subseteq A^-$. In particular, if $A^-$ is Lie solvable then so is $A^+$. Therefore, Theorem \ref{solvable} generalizes Theorem $1$ of \cite{Sal}. Moreover, since in characteristic $2$ the sets of symmetric and skew-symmetric elements obviously coincide, Examples 1 and 2 of \cite{Sal} show that the assumption $\char(\F)\neq 2$  in Theorems \ref{bengel}, $\ref{nilpotent}$,  and \ref{solvable}   cannot be removed.

\medskip
Finally, let $L$ be an arbitrary Lie algebra over a field $\F$ of characteristic $p\neq 2$ and denote by  $U(L)$ the ordinary universal enveloping algebra of $L$. Suppose that $U(L)^+$ (with respect to the principal involution) is Lie solvable or bounded Lie Engel. 
If $p=0$  then, by Theorem \ref{amitsur},  $U(L)$ satisfies a polynomial identity. Thus, $L$ is abelian    by a theorem of Latys\v{e}v (see e.g \cite[\S 6.7, Theorem 25]{Bah}).
Now suppose that  $p>2$ and consider the universal $p$-envelope 
of $L$
$$
\hat{L}:=\sum_{k \geq 0} {L^{p^k}}\subseteq U(L),
$$
where $L^{p^k}$ is the ${\mathbb F}$-vector space 
spanned by the set $\{l^{p^k}\vert\, l\in L \}$.
Then $\hat{L}$ is a restricted Lie algebra with $p$-map given by $h^{[p]}=h^p$ for all $h\in \hat{L}$. Note that by Corollary $1.1.4$ of \cite{S},  $U(L)=u(\hat{L})$. So,  Theorem \ref{solvable} and Theorem \ref{bengel} entail that $(\hat{L})^\prime$ is $p$-nilpotent. Since $u(\hat{L})=U(L)$ has no nontrivial zero divisors, we conclude that $L^\prime=0$. 

\smallskip
We have proved 
 
\begin{cor} Let $L$ be a Lie algebra over a field ${\mathbb F}$ of characteristic $p\neq 2$. Then $U(L)^+$ is Lie solvable or bounded Lie Engel if and only if $L$ is abelian.
\end{cor}  

\section*{acknowledgments}
The second author would like to thank the Department of  Matematics of the  University of  Salento for its hospitality during his visit while this work was completed.

\end{document}